\documentclass[12pt]{article}
\usepackage{amsmath,amsthm,amssymb,amscd}
\usepackage{latexsym}
\newtheorem{thm}{Theorem}[section]

 \newtheorem{prop}[thm]{Proposition}
 \theoremstyle{definition}
 
 \theoremstyle{remark}

 \numberwithin{equation}{section}

\title
{Greatest lower bounds on the transverse Ricci curvature of some toric Sasaki manifolds}

\author{ Hong Huang}
\date{}
\begin{document}
\maketitle
\begin{abstract}
 We determine the greatest lower bounds on the transverse Ricci curvature of compact toric Sasaki manifolds with positive basic first Chern class and with the first Chern class of the contact bundle being trivial. This is based on Wang-Zhu's and Futaki-Ono-Wang's works, and is an analogue of C. Li's work on toric Fano manifolds.



\end{abstract}
\maketitle


\section {Introduction}
In [S] Sz\'{e}kelyhidi defines the following invariant
\begin{equation*}
 R(X):=\sup \{t \hspace{1mm} | \hspace{1mm} \exists \hspace{1mm} \text {a} \hspace{1mm}  \text {K\"{a}hler} \hspace{1mm} \text {metric} \hspace{1mm} \omega \in c_1(X)\hspace{2mm} \text  {such} \hspace{1mm} \text  {that} \hspace{2mm} Ric(\omega)>t\omega\}
\end{equation*}
for a Fano manifold $X$.  (See also [T].) In [L] Li determines this  invariant for any compact toric Fano manifold $X$, based on Wang and Zhu's seminal work [WZ] on the existence of K\"{a}hler-Ricci soliton on any compact toric Fano manifold.
Note that recently Datar and Sz\'{e}kelyhidi [DS] recover the main results in [WZ] and [L] among other things, and Yao [Y] extends the result of [L] to the case of   homogeneous toric
bundles.

  In this note we first define an invariant analogous to $R(X)$ above for compact Sasaki manifolds with positive basic first Chern class and with the first Chern class of the contact bundle being trivial. Then, similarly to [L], we determine  the greatest lower bounds on the transverse Ricci curvature of compact toric Sasaki manifolds with positive basic first Chern class and with the first Chern class of the contact bundle being trivial, using and adapting Wang-Zhu's and Futaki-Ono-Wang's estimates in [WZ] and [FOW].

  As in for example [BGS] and [GZ], for a compact Sasaki manifold $S$ of dimension $2m+1$ with Sasaki structure $(\xi, \eta, \Phi, g)$ we define the space
\begin{equation*}
\mathcal{H}:=\{\phi \in C_B^\infty(S,\mathbb{R})  \hspace{1mm}|  \hspace{1mm} \eta_\phi=\eta+2d_B^c\phi \hspace{1mm}  \text {is}  \hspace{1mm} \text  {a}   \hspace{1mm} \text   {contact}    \hspace{1mm} \text  {form}\},
\end{equation*}
where $d_B^c = \frac{\sqrt{-1}}{2}(\bar{\partial}_B - \partial_B)$.   Assuming that   $c_1^B(S)>0$ and $c_1(D)=0$, where $c_1^B(S) $ is the basic first Chern class and $D=\text{Ker}$ $\eta$ is the contact bundle,
  following [S] we introduce the invariant
\begin{equation*}
\ \ \ R(S) := \sup \{t \hspace{1mm} | \hspace{1mm} \exists \hspace{1mm} \phi \in \mathcal{H}   \hspace{1mm}    \text {such} \hspace{1mm}  \hspace{1mm}    \text {that}   \hspace{1mm} \rho_\phi^T >t(m+1)d\eta_\phi \},
\end{equation*}
 where  $\rho_\phi^T$ is the transverse Ricci form derived from the Sasaki structure constructed in [FOW, Proposition 4.2] with transverse K\"{a}hler form $\frac{1}{2}d\eta_\phi$.

  Now we turn to the toric Sasaki manifolds; see for example [MS06], [MSY] and [FOW].
  Recall (see for example [FOW]) that a toric Sasaki manifold $S$ is a ($2m+1$)-dimensional Sasaki manifold with an effective action of a ($m+1$)-dimensional  torus $G\cong T^{m+1}$ preserving the Sasaki structure $(\xi,\eta,\Phi,g)$ such that the Reeb field $\xi$ is induced by an element of the Lie algebra $\mathfrak{g}$ of $G$. Thus the   cone $(C(S), \bar{g})=(\mathbb{R}_+\times S, dr^2+r^2g)$ of a toric Sasaki manifold $S$ is a toric K\"{a}hler manifold.

  Let $S$ be a ($2m+1$)-dimensional  compact toric Sasaki manifold. The moment map $\mu_\eta:S\rightarrow \mathfrak{g}^*$ w.r.t. the contact form $\eta$ is given by
  \begin{equation*}
  \langle \mu_\eta(x), X\rangle=\eta(X_S(x)), \hspace{4mm}  \forall x\in S,
  \end{equation*}
  where $X_S$ is the vector field on $S$ induced by $X\in \mathfrak{g}$, i.e., $X_S(x):= \frac{d}{dt}|_{t=0}exp(tX)\cdot x$. On the other hand, the complexification $G^c\cong (\mathbb{C}^*)^{m+1}$ acts on $C(S)$ by biholomorphic automorphisms, and the corresponding  moment map  $\mu:C(S)\rightarrow \mathfrak{g}^*$  w.r.t. the  K\"{a}hler form $\omega=d(\frac{1}{2}r^2\eta)(=dr^2+r^2g)$
  (here the pull-back of $\eta$ by the projection  $C(S)\rightarrow S$ is still denoted by $\eta$) is given by
  \begin{equation*}
  \langle \mu(x), X\rangle=r^2\eta(X_S(x)), \hspace{4mm} \forall x\in C(S),
  \end{equation*}
   where $X_S$ is  viewed as a vector field on $C(S)$. We denote the image of $\mu$ by $C(\mu)$, which is a convex rational polyhedral cone.  So there exist vectors $\lambda_a, a=1,\cdot\cdot\cdot,d,$ in the integral lattice $\mathbb{Z}_\mathfrak{g}:=$Ker$\{\exp: \mathfrak{g} \rightarrow G\}$ such that
   \begin{equation*}
   C(\mu)=\{y\in \mathfrak{g}^* \hspace{2mm} |\hspace{2mm} l_a(y)=\langle y, \lambda_a\rangle \geq 0, \hspace{2mm} a=1,\cdot\cdot\cdot,d\}.
   \end{equation*}

   We also denote  the interior of $C(\mu)$ by Int$C(\mu)$. It is easy to see that the image of $\mu_\eta$
   \begin{equation*}
   \text{Im}(\mu_\eta)=\{\alpha \in C(\mu)\hspace{1mm} |\hspace{1mm} \alpha(\xi)=1\}.
   \end{equation*}

Now we assume further that the compact toric Sasaki manifold $S$ has  $c_1^B(S)>0$ and $c_1(D)=0$.  Then by [FOW, Proposition 4.3],  $c_1^B(S)$ is represented by $\tau d\eta$ for some positive constant $\tau$. Using  $\mathcal{D}$-homothetic transformation if needed we may and will assume that  $(m+1)d\eta \in 2\pi c_1^B(S)$. Moreover, by [FOW, Proposition 6.7], using transverse K\"{a}hler deformation if needed we may and will further assume that the symplectic potential on $(C(S), d(\frac{1}{2}r^2\eta))$ is given by formula (42) in [FOW]. Then by [FOW] there exists a unique rational vector $\gamma \in \mathfrak{g}^*$ such that
\begin{equation*}
\langle \gamma, \lambda_a\rangle=-1,  \hspace{4mm}  a=1,\cdot\cdot\cdot,d.
 \end{equation*}

  Choose a $m$-dimensional subtorus $H\subset G$ whose Lie algebra is
  \begin{equation*}
  \mathfrak{h}:=\{x \in \mathfrak{g}\hspace{1mm} | \hspace{1mm} \langle \gamma, x\rangle=0\}.
   \end{equation*}
    In particular, since $\langle \gamma, \xi\rangle=-(m+1)$ ([FOW, (49)]),  $\mathfrak{h}$ does not contain  $\xi$. (Here we have identified the Reeb field $\xi$ with the element in $\mathfrak{g}$ which induces it.) Let $H^c\cong (\mathbb{C}^*)^m$ be the complexification of $H$. Fix a point $p \in \mu^{-1}$(Int$C(\mu)$), let $Orb_{C(S)}(H^c,p)$ be the orbit through $p$ of the $H^c$-action on $C(S)$. The  moment map $\mu_{\eta,H}: Orb_{C(S)}(H^c,p) \rightarrow \mathfrak{h}^*$ on the K\"{a}hler manifold $(Orb_{C(S)}(H^c,p), \frac{1}{2}d\eta|_{Orb_{C(S)}(H^c,p)})$  for the $H$-action is defined  by
  \begin{equation*}
  \langle\mu_{\eta,H}(y), X\rangle=\eta(X)(y),   \hspace{4mm} \hspace{2mm} y \in Orb_{C(S)}(H^c,p), X \in \mathfrak{h},
  \end{equation*}
  where the $X$ on the RHS of the equality is the vector field on $Orb_{C(S)}(H^c,p)$ induced by $X \in \mathfrak{h}$.  It turns out that
  \begin{equation*}
  \text{Im}(\mu_{\eta,H})=\iota^*(\text{Im}(\mu_\eta))=\{\iota^*\alpha \hspace{2mm} | \hspace{2mm}  \alpha \in C(\mu), \alpha(\xi)=1\},
  \end{equation*}
  where $\iota: \mathfrak{h} \rightarrow \mathfrak{g}$ is the inclusion map. (See [FOW].) This image is a compact convex polyhedron.
  It is not necessarily rational, since the Sasaki structure on $S$ may not be quasi-regular. (Compare [MS06] and [FOW].)

   On $Orb_{C(S)}(H^c,p)\cong (\mathbb{C}^*)^m$ we introduce the affine logarithm coordinates
  \begin{equation*}
  (w^1,\cdot\cdot\cdot,w^m)=(x^1+\sqrt{-1}\theta^1,\cdot\cdot\cdot,x^m+\sqrt{-1}\theta^m)
  \end{equation*}
  for a point
   \begin{equation*}
   (e^{x^1+\sqrt{-1}\theta^1},\cdot\cdot\cdot,e^{x^m+\sqrt{-1}\theta^m})\in (\mathbb{C}^*)^m \cong Orb_{C(S)}(H^c,p).
   \end{equation*}
  Now $\frac{1}{2}d\eta|_{Orb_{C(S)}(H^c,p)}$ is determined by a convex function $u^0$ on $\mathbb{R}^m$,
  \begin{equation*}
  \frac{1}{2}d\eta|_{Orb_{C(S)}(H^c,p)}=\sqrt{-1}\partial\bar{\partial}u^0=\frac{\sqrt{-1}}{4}\frac{\partial^2u^0}{\partial x^i\partial x^j}dw^i\wedge d\overline{w^j}.
  \end{equation*}
  It is easy to see (cf. for example [FOW]) that  (after translation)  the interior Int(Im($\mu_{\eta,H}$))  can be identified with
  \begin{equation*}
  \Sigma:=\{Du^0(x)=(\frac{\partial u^0}{\partial x^1}(x),\cdot\cdot\cdot,\frac{\partial u^0}{\partial x^m}(x))\hspace{2mm} | \hspace{2mm}  x\in \mathbb{R}^m\}.
\end{equation*}
We call the closure  $\overline{\Sigma}$   the moment polytope  of $(Orb_{C(S)}(H^c,p), \frac{1}{2}d\eta|_{Orb_{C(S)}(H^c,p)})$ for the $H$-action (compare (57) in [FOW]).

 It follows from  Proposition  7.3 (or Lemma 7.5) of [FOW] that the origin $O$ of $\mathbb{R}^m$ is contained in $\Sigma$. We observe that the  barycenter $P_c$ of the moment polytope $\overline{\Sigma}$ coincides with the origin $O$ if and only if  the Sasaki-Futaki invariant $f$ of $S$ (for definition see [BGS] and [FOW]) vanishes, see Proposition 3.4.

Similarly to [L, Theorem 1] we have
\begin{thm} \label{thm 1.1} \ \
 Let $(S, \xi, \eta,\Phi,g)$ be a compact ($2m+1$)-dimensional toric Sasaki manifold with positive basic first Chern class and with the first Chern class of the contact bundle being trivial.
  Let $\overline{\Sigma}$ and $P_c$ be as above.

  \noindent If $P_c\neq O$,  then $R(S)<1$,  and
 \begin{equation*}
 R(S)=\frac{|\overline{OQ}|}{|\overline{P_cQ}|},
 \end{equation*}
where $Q$ is the intersection of the ray $P_c+\mathbb{R}_{\geq 0}\cdot \overrightarrow{P_cO}$ with $\partial \overline{\Sigma}$.

\noindent If $P_c =O$, then $S$ admits a Sasaki-Einstein metric, and $R(S)=1$.
\end{thm}

That if $P_c =O$ then $S$ admits a Sasaki-Einstein metric  follows from [FOW], see also the proof of Proposition 3.4 below;  we include it here for completeness. The bridge between $R(S)$ and $\frac{|\overline{OQ}|}{|\overline{P_cQ}|}$ is Aubin's continuity path for finding Sasaki-Einstein metrics.
In Section 2, following [S] we show that on a ($2m+1$)-dimensional compact Sasaki manifold $(S, \xi)$ (not necessarily toric) with positive basic first Chern class and with the first Chern class of the contact bundle being trivial, $R(S)$ is equal to the maximum existence time of Aubin's continuity path for  finding Sasaki-Einstein metrics on $S$. In Section 3 we use this continuity path to prove Theorem 1.1.

For the most part of the proof of Theorem 1.1 we follow closely the lines of [L] (see also [Y]), using and/or adapting estimates from [FOW] and [WZ]. However, there is one point where our argument is  slightly different from that in [L]: To prove the Claim 1 on p.4929 of [L], Li uses the simple formula (2) on p. 4923 of [L] expressing  the initial K\"{a}hler potential $\tilde{u}_0$ via the vertices of the moment polytope. In our case, such a simple expression for $u^0$ is not available in general (when the Sasaki structure is not quasi-regular). Instead we have the formula (81) on p. 621 of [FOW] for $u^0$, which is somewhat difficult to treat  directly.  The idea is to use the Legendre transform to convert the K\"{a}hler potential $u^0$ to the symplectic potential $G_0(v)$,  and exploit the degenerate behavior of  (Hess $G_0(v))^{-1}$ near the boundary  $\partial \overline{\Sigma}$ to  prove a result similar to the Claim 1 in [L]. (Compare also [FOW] and [D].)

\section{The invariant $R(S)$}

Let $(S, \xi, \eta, \Phi, g)$  be  a ($2m+1$)-dimensional compact (not necessarily toric) Sasaki manifold  with positive basic first Chern class, with $c_1(D)=0$ ($D= \text {Ker}$ $\eta$) and with $(m+1)d\eta \in 2\pi c_1^B(S)$.

Define (cf. for example [FOW], [GZ], [Z11a]) Mabuchi functional on $\mathcal{H}$ (see the Introduction) via its variation
\begin{equation*}
\frac{d}{dt}\mathcal{M}(\phi_t)=\int_S\dot{\phi_t}(2m(m+1)-s_{\phi_t}^T)(\frac{1}{2}d\eta_{\phi_t})^m\wedge \eta
\end{equation*}
and the requirement $\mathcal{M}(0)=0$, where $s_{\phi_t}^T$ is the transverse scalar curvature derived from the Sasaki structure constructed in  [FOW, Proposition 4.2] with transverse K\"{a}hler form $\frac{1}{2}d\eta_{\phi_t}$ .

Let $\chi$ be a transverse K\"{a}hler form on $S$. We also define (cf. for example [VZ]) the $\mathcal{J}_\chi$  functional on  $\mathcal{H}$   via its variation
\begin{equation*}
\frac{d}{dt}\mathcal{J}_\chi(\phi_t)=2m(m+1)\int_S\dot{\phi_t}( \chi \wedge(\frac{1}{2}d\eta_{\phi_t})^{m-1}-(\frac{1}{2}d\eta_{\phi_t})^m)\wedge \eta
\end{equation*}
and the requirement $\mathcal{J}_\chi(0)=0$.  Compare also [S].

Given $\psi \in \mathcal{H}$, let $h_\psi$ be determined by
\begin{equation*}
\rho_\psi^T-(m+1)d\eta_\psi=\sqrt{-1}\partial_B\bar{\partial}_Bh_\psi
\end{equation*}
and
\begin{equation*}
\int_S e^{h_\psi}(\frac{1}{2}d\eta_\psi)^m\wedge \eta=\int_S (\frac{1}{2}d\eta_\psi)^m\wedge \eta,
\end{equation*}
Aubin's continuity path for finding Sasaki-Einstein metrics is given by the following transverse Monge-Amp\`{e}re equation for $\phi_t \in \mathcal{H}$
\begin{equation*}
\frac{(d\eta+2\sqrt{-1}\partial_B\bar{\partial}_B\phi_t)^m\wedge \eta}{(d\eta+2\sqrt{-1}\partial_B\bar{\partial}_B\psi)^m\wedge \eta}=e^{h_\psi-t(2m+2)\phi_t},
\end{equation*}
or
\begin{equation*}
 \frac{\det (g_{i\bar{j}}^T+\phi_{i\bar{j}}) }  {\det (g_{i\bar{j}}^T+\psi_{i\bar{j}})}=\exp (h_\psi-t(2m+2)\phi_t).   \hspace{12mm}   (*)_t
\end{equation*}
The equation $(*)_t$ is  equivalent to $\rho_{\phi_t}^T=(m+1)(td\eta_{\phi_t}+(1-t)d\eta_\psi)$.  When $t=0$ the equation is solvable by the transverse Yau theorem in [E].

Following [S] we call a functional $\mathcal{F}$ defined on the space $\mathcal{H}$ proper if there exist constants $\epsilon, C>0$ such that
\begin{equation*}
\mathcal{F}(\psi)> \epsilon \mathcal{J}_{\frac{1}{2}d\eta}(\psi)-C
\end{equation*}
for any $\psi \in \mathcal{H}$.

\begin{thm} \label{thm 2.1} \ \
 Let $S$ be as above. The following are equivalent for $0\leq t<1$.

 1) Given any $\psi\in \mathcal{H}$  the equation $(*)_t$ can be solved.

 2) There exists $\psi \in \mathcal{H}$ such that $\hspace{1mm} \rho_\psi^T >t(m+1)d\eta_\psi$.

 3) The functional $\mathcal{M}+(1-t)\mathcal{J}_{\frac{1}{2}d\eta_\psi}$ is proper for any $\psi \in \mathcal{H}$.
\end{thm}

\noindent {\bf Proof} \ \ The proof is along the lines of proof of Theorem 1 in [S].  We only indicate some necessary modifications.  We use [JZ] and [vC] to replace  [CT08] in the proof of Proposition 3 in [S], and
 use [NS] to replace [BM87] in the proof of Lemma 5 in [S].   \hfill{$\Box$}

\section{ Proof of Theorem 1.1}

Let $(S, \xi, \eta,\Phi,g)$  be a compact ($2m+1$)-dimensional toric Sasaki manifold satisfying the assumptions of Theorem 1.1.
  Choose $H$, $p$ and $u^0$ as in the Introduction.

Choose $\psi=0$ in $(*)_t$ of Section 2. Then as in [FOW], $(*)_t$ can be converted to  the following Monge-Amp\`{e}re equation for a strictly convex function $u$
\begin{equation*}
\det (u_{ij})=\exp (-(2m+2)(tu+(1-t)u^0))        \hspace{4mm}     \text { on } \hspace{2mm} \mathbb{R}^m.              \hspace{8mm} (**)_t
\end{equation*}
By Theorem 2.1, $(**)_t$ is solvable when $t< R(S)$.

Let $u$ be a solution to $(**)_t$, and
\begin{equation*}
w_t=tu+(1-t)u^0.
\end{equation*}

\noindent Since $Dw_t(\mathbb{R}^m)=Du(\mathbb{R}^m)=Du^0(\mathbb{R}^m)=\Sigma $ (compare for example, [M1], [WZ], and the proof of Fact 2 in Section 2 of [Hu1]) and $O\in \Sigma $ (as observed in the Introduction), the strictly convex function $w_t$ is proper, and  attains its  minimum $m_t$ at a unique point $x_t \in \mathbb{R}^m$.

\begin{prop} \label{prop 3.1} \  \  1) There exists a constant $C$ independent of $t < R(S)$, such that
\begin{equation*}
|m_t|\leq C.
\end{equation*}
2) There exist $\kappa >0$ and a constant $C$, both independent of $t < R(S)$, such that
\begin{equation*}
w_t\geq \kappa |x-x_t|-C.
\end{equation*}
\end{prop}

\noindent {\bf Proof} The proof is the same as that of Proposition 2 in [L], which uses arguments of [WZ, Lemma 3.2] and [D, Section 3.4, Proposition 1]; compare also the proof of Lemma 3.1 in [Hu1].
\hfill{$\Box$}

\vspace*{0.4cm}

\begin{prop} \label{prop 3.2} \  \
Fix $t_0$. There exists a constant $C_1$ such that $|x_t|\leq C_1$  for $0\leq t \leq t_0$, where $x_t$ is the minimum point of $w_t=tu+(1-t)u^0$ with $u$ being any solution to $(**)_t$ if and only if there exists a constant $C_2$ such that $|\varphi_t|\leq C_2$ for $0\leq t \leq t_0$, where $u^0+\varphi_t$ is any solution to $(**)_t$.
\end{prop}
\noindent {\bf Proof} \ \  The proof is similar to that of [L, Proposition 3] with the help of Proposition 3.1, 1),
[FOW, Proposition 7.3] and [Z11b, Theorem 1.1].  (Alternatively, one can also use Proposition 3.1, 2) and the argument in the last paragraph of Section 3 in [D].)
  \hfill{$\Box$}

\vspace*{0.4cm}

\begin{prop} \label{prop 3.3}
If $R(S)<1$,  there exist a sequence $\{t_k\}$ and a point $y_\infty \in \partial\overline{\Sigma}$, such that
\begin{equation*}
\lim_{k\rightarrow \infty} t_k=R(S),    \hspace{4mm}  \lim_{k\rightarrow \infty} |x_{t_k}|=\infty,   \hspace{4mm}  \lim_{k\rightarrow \infty}Du^0(x_{t_k})=y_\infty.
\end{equation*}
\end{prop}
\noindent {\bf Proof} \ \  The result follows easily from Theorem 2.1, Proposition 3.2, the properness of $u^0$ and the compactness of $\overline{\Sigma}$.
  \hfill{$\Box$}

\vspace*{0.4cm}

Recall [FOW] that
\begin{equation}
\overline{\Sigma}=\cap_{a=1}^d\{l_a'(v)\geq0\},
\end{equation}
where $l_a'(v)=\langle v, \lambda_a'\rangle +\frac{1}{m+1}$,  and $\lambda_a' \in \mathfrak{h}\cong \mathbb{R}^m$ is given by the decomposition
\begin{equation*}
\lambda_a=\iota (\lambda_a')+\frac{1}{m+1}\xi,
\end{equation*}
where $\lambda_a$ is as in the Introduction.

W.l.o.g. we may assume that

$\begin {array}{l}
l_a'(y_\infty)=0, \hspace{4mm}   a=1, \cdot\cdot\cdot, d_0, \\
l_a'(y_\infty)>0,  \hspace{4mm}  a=d_0+1, \cdot\cdot\cdot, d,
\end {array} $

\noindent where $ d_0 \geq 1$.

Note that  we have
\begin{equation}
\int_{\mathbb{R}^m} e^{-(2m+2)w_t}dx=\int_{\mathbb{R}^m} \det(u_{ij})dx=\int_{\Sigma}dy=Vol(\Sigma).
\end{equation}

Since $w_t$ is a proper strictly convex function on ${\mathbb{R}^m}$,  $w_t(x)\rightarrow +\infty$ as $|x|\rightarrow \infty$.
So we have
\begin{equation*}
\int_{\mathbb{R}^m} \frac{\partial w_t}{\partial x^i}e^{-(2m+2)w_t}dx=-\frac{1}{2m+2}\int_{\mathbb{R}^m} \frac{\partial e^{-(2m+2)w_t}}{\partial x^i}dx=0,\hspace{4mm}   i=1,\cdot\cdot\cdot,m.
\end{equation*}
(Compare also for example [D].)
It follows that when $t<R(S)$,
\begin{equation*}
\int_{\mathbb{R}^m} (Du^0)e^{-(2m+2)w_t}dx=-\frac{t}{1-t}\int_{\mathbb{R}^m} (Du)e^{-(2m+2)w_t}dx.
\end{equation*}
On the other hand,
\begin{equation*}
\int_{\mathbb{R}^m} (Du)e^{-(2m+2)w_t}dx=\int_{\mathbb{R}^m} (Du)\det(u_{ij})dx=\int_{\Sigma}ydy=Vol(\Sigma)P_c,
\end{equation*}
where $P_c$ is the barycenter of $\overline{\Sigma}$ (as in the statement of Theorem 1.1).

Thus as in [L] we get
\begin{equation}
\frac{1}{Vol(\Sigma) }\int_{\mathbb{R}^m} (Du^0)e^{-(2m+2)w_t}dx=-\frac{t}{1-t}P_c
\end{equation}
when $t<R(S)$.

As in [Y] we define
\begin{equation}
R_\Sigma:=\sup\{t \hspace{2mm}| \hspace{2mm} 0\leq t <1, -\frac{t}{1-t}P_c \in \Sigma \}.
\end{equation}
Since   $e^{-(2m+2)w_t}>0$, $\frac{1}{Vol(\Sigma) }\int_{\mathbb{R}^m} e^{-(2m+2)w_t}dx=1$ by (3.2), $Du^0(x) \in \Sigma$ for any $x \in \mathbb{R}^m$, and $\Sigma$ is convex, we have
\begin{equation*}
\frac{1}{Vol(\Sigma) }\int_{\mathbb{R}^m} (Du^0)e^{-(2m+2)w_t}dx \in \Sigma.
\end{equation*}
 Combining with (3.3) we get that $-\frac{t}{1-t}P_c \in \Sigma$ for $t < R(S)$.  So
 \begin{equation}
 R(S) \leq R_\Sigma.
 \end{equation}
 Compare [Y]. In particular, if $R(S)=1$, then $R_\Sigma=1$.

 Let $(w^1,\cdot\cdot\cdot,w^m)$ be  the affine logarithm coordinates on $Orb_{C(S)}(H^c,p)\cong (\mathbb{C}^*)^m$ as in the introduction, let $X_k=-\frac{\sqrt{-1}}{2}\frac{\partial}{\partial w^k}$,  and  $\theta_{X_k}$ be its Hamiltonian function  (see p.597 and p.604 of [FOW]), $k=1,\cdot\cdot\cdot,m$. By [FOW, Lemma 7.4], $\theta_{X_k}=\frac{\partial u^0}{\partial x^k}$.

The following result is implicitly from [FOW], and is analogous to  [M1, Corollary  5.5], [M2, Lemma 6.1], and [F, Theorem 3.4.1].
\begin{prop} \label{prop 3.4} \  \ (cf. [FOW])
  $P_c=O$ if and only if  the  Sasaki-Futaki invariant of $S$  vanishes.
\end{prop}
\noindent {\bf Proof} \ \ For $k=1,\cdot\cdot\cdot,m$,  we compute as in the proof of [FOW, Lemma 7.5],

$$\begin {array}{l}
f(X_k)=-\int_S \theta_{X_k}(\frac{1}{2}d\eta)^m\wedge \eta \\
=-\int_S \frac{\partial u^0}{\partial x^k} \det (u_{ij}^0) dx \wedge d\theta \wedge \eta \\
=-\text {const.} \int_\Sigma y_kdy.
\end {array} $$

\noindent So that the Sasaki-Futaki invariant $f$ of $S$  vanishes implies that $\int_\Sigma y_kdy=0$ for $1\leq k \leq m$, and $P_c=O$.

On the other hand, if  $P_c=O$,  then $c_i=0$ ($1\leq i \leq m$) satisfy the equations in Lemma 7.5 in [FOW]. Since that the $c_i$ ($1\leq i \leq m$) satisfying the equations in Lemma 7.5 in [FOW] are unique (compare [TZ, Lemma 2.2] and the Remark on p. 93 of [WZ]), we see that
the vector field $X$ in Proposition 5.3 of [FOW] must be trivial. Then from the proof of  [FOW,Theorem 1.1] we see that $S$ admits a Sasaki-Einstein metric, and the Sasaki-Futaki invariant of $S$ vanishes.    \hfill{$\Box$}

\vspace*{0.4cm}

\begin{prop} \label{prop 3.5} \  \
Suppose   $R(S)< 1$.  Let  $Q:=-\frac{R(S)}{1-R(S)}P_c$, then $Q \in \partial \overline{\Sigma}$. More precisely $Q$ lies on the same faces of $\overline{\Sigma}$ as the point $y_\infty$ does, that is,

$$\begin {array}{l}
l_a'(Q)=0, \hspace{4mm}   a=1, \cdot\cdot\cdot, d_0, \\
l_a'(Q)>0,  \hspace{4mm}  a=d_0+1, \cdot\cdot\cdot, d.
\end {array} $$

\noindent Consequently in this case $P_c\neq O$.
\end{prop}

\noindent {\bf Proof}. \ \
Using (3.3), (3.2) and (3.1) we get

$$\begin {array}{l}
l_a'(-\frac{t}{1-t}P_c)=\frac{1}{vol(\Sigma)}\int_{\mathbb{R}^m} \langle Du^0, \lambda_a' \rangle e^{-(2m+2)w_t}dx+\frac{1}{m+1}\\
=\frac{1}{vol(\Sigma)}\int_{\mathbb{R}^m} (\langle Du^0, \lambda_a' \rangle +\frac{1}{m+1})e^{-(2m+2)w_t}dx \geq 0.
\end {array} $$

Since $R(S)<1$, we can let $t \rightarrow R(S)$ and get that
\begin{equation*}
l_a'(-\frac{R(S)}{1-R(S)}P_c)\geq 0,  \hspace{4mm}     a=1,\cdot\cdot\cdot,d.
\end{equation*}

Now the rest of the arguments is almost the same as in the proof of Proposition 4 in [L], using Propositions 3.1,  2),  and Proposition 3.3, with the  Claim 1 in [L]  replaced by the Claim below.

\vspace*{0.4cm}

\noindent {\bf Claim }  (Compare p. 56 of [D]) The derivative of the function $s_a(x):=\log (l_a'(Du^0(x)))$ is bounded on $\mathbb{R}^m$.

\noindent {\bf Proof of Claim.} \ \ We compute
\begin{equation*}
Ds_a(x)=\frac{D^2u^0(x)\lambda_a'}{l_a'(Du^0(x))}=\frac{(D^2G_0(v))^{-1}\lambda_a'}{l_a'(v)},
\end{equation*}
where $v=Du^0(x)$, and $G_0(v)$ is the Legendre transform of the K\"{a}hler potential $u^0(x)$, and is the symplectic potential of the  the K\"{a}hler manifold

$(Orb_{C(S)}(H^c,p), \frac{1}{2}d\eta|_{Orb_{C(S)}(H^c,p)})$.

By computing the Hessian $D^2G_0(v)$ using  formula (82) of [FOW] one sees that  as one approaches the $(m-1)$-dimensional face $l_a'(v)=0$ of $\overline{\Sigma}$ from the interior, the positive definite matrix $(D^2G_0(v))^{-1}$ will tend to be degenerate, and will acquire a kernel that is generated by the normal $\lambda_a'$ when one reaches the face $l_a'(v)=0$ at last.  (Compare for example [A] and the proof of Fact 3 in Section 2 of [Hu1].)
  So $\frac{(D^2G_0(v))^{-1}\lambda_a'}{l_a'(v)}$ can be extended to a  continuous  function on the closure $\overline{\Sigma}$.    \hfill{$\Box$}

\vspace*{0.4cm}

Combining (3.5) with Proposition 3.5 we get that $R(S)=R_\Sigma$. Now we see that if $P_c \neq O$, then $R_\Sigma < 1$ by definition of $R_\Sigma$ (see (3.4)) and the compactness of $\overline{\Sigma}$, and $R(S) < 1$.
If $P_c = O$, then $R_\Sigma = 1$  by definition and the fact $O \in \Sigma$, and $R(S)=1$.

Now Theorem 1.1 is proved.

\vspace*{0.4cm}

{\bf Remark 1}.   The statement that $R(S) < 1$  implies that $P_c\neq O$ can also be proved as follows: If $P_c= O$, then by [FOW] (see the proof of Proposition 3.4 here) there is a Sasaki-Einstein metric on $S$, which implies
$R(S) = 1$. So that $R(S) < 1$  implies that $P_c\neq O$.

\vspace*{0.4cm}

{\bf Remark 2}.  In the proof of Proposition 4 in [L], Li uses the fact that in his situation, $P_c\neq O$ implies $R(X_\triangle)<1$, although he does not state it explicitly. This fact can be easily deduced from the arguments in [L]; compare  [Y] (and the above proof of Theorem 1.1). It can also be proved as follows: If $P_c\neq O$, then by [M1] the Futaki invariant of $X_\triangle$ does not vanish, and $X_\triangle$ cannot be $K$-semistable, so $R(X_\triangle) <1$ by Corollary 1.1 of [MS].

In our situation, the fact that if $P_c\neq O$ then $R(S)<1$ (which was proved above) can also be proved as follows: By Proposition 3.4, the Sasaki-Futaki invariant of $S$ does not vanish.  So by [BHLT] $S$ can not be $K$-semistable. In a forthcoming paper [Hu2] we'll show
that for  a ($2m+1$)-dimensional compact (not necessarily toric) Sasaki manifold $S$  with positive basic first Chern class and  with $c_1(D)=0$ ($D= \text {Ker}$ $\eta$), $R(S)=1$ implies that $S$ is $K$-semistable. (The proof is along the lines of [MS], uses Sasaki-Ricci flow (cf. [C], [H]), and also uses  [CS], [JZ] and [vC].)
So that $P_c\neq O$  implies that $R(S)<1$.

\vspace*{0.4cm}

\noindent {\bf Acknowledgements} {\hspace*{4mm}}     I was partially supported by NSFC
(No. 11171025) and  Laboratory of Mathematics and Complex Systems, Ministry of Education.


\hspace *{0.4cm}

\bibliographystyle{amsplain}

\noindent {\bf References}

\hspace *{0.1cm}

\vspace*{0.4cm}

[A] M. Abreu, K\"{a}hler geometry of toric manifolds in symplectic coordinates,  Symplectic and contact topology: interactions and perspectives,  1-24, Fields Inst. Commun., 35, Amer. Math. Soc., Providence, RI, 2003.

[BM87] S. Bando and T. Mabuchi, Uniqueness of Einstein K\"{a}hler metrics modulo connected group actions, in Algebraic geometry, Sendai, 1985, Advanced Studies in Pure Mathematics, vol. 10 (North-Holland, Amsterdam, 1987), 11-40.

[BHLT] C. Boyer, H.N. Huang, E. Legendre, C. T{\o}nnesen-Friedman, The Einstein-Hilbert functional and the Sasaki-Futaki invariant,   Int. Math. Res. Not., 2017, no.7,  1942-1974.

[BGS]  C. Boyer, K. Galicki, S. Simanca, Canonical Sasakian metrics. Comm. Math. Phys.  279  (2008),  no. 3, 705-733.

 [CT08] X. X. Chen, G. Tian, Geometry of K\"{a}hler metrics and foliations by holomorphic discs, Publ. Math. Inst. Hautes \'{E}tudes Sci. 107 (2008), 1-107. 

[C]  T. Collins,  The transverse entropy functional and the Sasaki-Ricci flow. Trans. Amer. Math. Soc.  365  (2013),  no. 3, 1277-1303.

[CS] T. Collins,  G. Sz\'{e}kelyhidi, $K$-Semistability for irregular Sasakian manifolds,   arXiv:1204.2230, to appear in J. Diff. Geom.

[DS] V. Datar, G. Sz\'{e}kelyhidi, K\"{a}hler-Einstein metrics along the smooth continuity method,  Geom. Funct. Anal. 26 (2016), no. 4, 975-1010.

[D] S. Donaldson, K\"{a}hler geometry on toric manifolds, and some other manifolds with large symmetry.  Handbook of geometric analysis. No. 1,  29-75, Adv. Lect. Math. (ALM), 7, Int. Press, Somerville, MA, 2008.

[E] A. El Kacimi-Alaoui, Op\'{e}rateurs transversalement elliptiques sur un feuilletage riemannien et applications. (French)  Compositio Math.  73  (1990),  no. 1, 57-106.

[F] A. Futaki, K\"{a}hler-Einstein metrics and integral invariants.
Lecture Notes in Mathematics, 1314. Springer-Verlag, Berlin, 1988.

[FOW]  A. Futaki, H. Ono, G. Wang, Transverse K\"{a}hler geometry of Sasaki manifolds and toric Sasaki-Einstein manifolds.
J. Diff. Geom.  83  (2009),  no. 3, 585-635.

[GZ]  P. Guan, X. Zhang,  Regularity of the geodesic equation in the space of Sasakian metrics,  Adv. Math.  230  (2012),  no. 1, 321-371.

[H] W. He, The Sasaki-Ricci flow and compact Sasaki manifolds of positive transverse holomorphic bisectional curvature. J. Geom. Anal.  23  (2013),  no. 4, 1876-1931.

[Hu1] H. Huang, K\"{a}hler-Ricci flow on homogeneous toric bundles, arXiv:1705.07735.

[Hu2] H. Huang, Sasaki-Ricci flow and $K$-semistability, in preparation.

[JZ] X. Jin, X. Zhang, Uniqueness of constant scalar curvature Sasakian metrics,  Ann. Glob. Anal. Geom.  49 (2016), 309-328.

[L] C. Li, Greatest lower bounds on Ricci curvature for toric Fano manifolds. Adv. Math.  226  (2011),  no. 6, 4921-4932.

[M1] T. Mabuchi, Einstein-K\"{a}hler forms, Futaki invariants and convex geometry on toric Fano varieties. Osaka J. Math.  24  (1987),  no. 4, 705-737.

 [M2] T. Mabuchi,  An algebraic character associated with the Poisson brackets. Recent topics in differential and analytic geometry, 339-358, Adv. Stud. Pure Math., 18-I, Academic Press, Boston, MA, 1990.

[MS06] D. Martelli, J. Sparks, Toric geometry, Sasaki-Einstein manifolds and a new infinite class of AdS/CFT duals. Comm. Math. Phys.  262  (2006),  no. 1, 51-89.

[MSY] D. Martelli, J. Sparks, S. T. Yau,
 The geometric dual of a -maximisation for toric Sasaki-Einstein manifolds.
Comm. Math. Phys.  268  (2006),  no. 1, 39-65.

 [MS] O. Munteanu,  G. Sz\'{e}kelyhidi, On convergence of the K\"{a}hler-Ricci flow,  Comm. Anal. Geom. 19 (2011),  887-903.

[NS] Y. Nitta, K. Sekiya, Uniqueness of Sasaki-Einstein metrics,   Tohoku Math. J. (2)  64 (2012), no.3, 453-468.

[S] G. Sz\'{e}kelyhidi, Greatest lower bounds on the Ricci curvature of Fano manifolds,
Compos. Math.  147  (2011),  no. 1, 319-331.

[T] G. Tian, On stability of the tangent bundles of Fano varieties, Internat. J. Math. 3 (1992), no.3, 401-413.

[TZ] G. Tian, X. H. Zhu, A new holomorphic invariant and uniqueness of K\"{a}hler-Ricci solitons, Comment. Math. Helv. 77 (2002), 297-325.

[vC] C. van Coevering, Monge-Amp\`{e}re operators, energy functionals, and uniqueness of Sasaki-extremal metrics, arXiv:1511.09167.

[VZ] L. Vezzoni, M. Zedda,  On the J-flow in Sasakian manifolds,  Ann. di Mat. Pura e Appl. 195 (2016), 757-774.

[WZ] X. Wang, X. H. Zhu, K\"{a}hler-Ricci solitons on toric manifolds with positive first Chern class.
Adv. Math.  188  (2004),  no. 1, 87-103.

[Y] Y. Yao, Greatest lower bounds on Ricci curvature of homogeneous toric bundles, Internat. J.  Math. 28 (2017), no. 4, 1750024 (16 pages).

[Z11a]  X. Zhang, Energy properness and Sasakian-Einstein metrics. Comm. Math. Phys.  306  (2011),  no. 1, 229-260.

[Z11b] X. Zhang,  Some invariants in Sasakian geometry. Int. Math. Res. Not.  2011,  no. 15, 3335-3367.

\vspace *{0.4cm}

School of Mathematical Sciences, Beijing Normal University,

Laboratory of Mathematics and Complex Systems, Ministry of Education,

Beijing 100875, P.R. China

 E-mail address: hhuang@bnu.edu.cn

\end{document}